# Simultaneous Sizing of a Rocket Family with Embedded Trajectory Optimization

Byeong-Un Jo[*]
*Sejong University, Seoul, Republic of Korea, 05006*
and
Koki Ho[†]
*Georgia Institute of Technology, Atlanta, Georgia, 30332*

**This paper presents a sizing procedure for a rocket family capable of fulfilling multiple missions considering the commonalities between the vehicles. The procedure aims to take full advantage of sharing a common part across multiple rockets whose payload capability differs entirely, ultimately leading to cost savings in designing a rocket family. As the foundation of the proposed rocket family design method, an integrated sizing method with trajectory optimization for a single rocket is first formulated as a single optimal control problem. This formulation can find the optimal sizing along with trajectory results in a tractable manner. Building upon this formulation, the proposed rocket family design method is developed to (1) determine the feasible design space of the rocket family design problem (i.e., commonality check); and (2) if a feasible design space is determined to exist, minimize the cost function within that feasible space by solving an optimization problem in which the optimal control problem is embedded as a subproblem. A case study is carried out on a rocket family composed of expendable and reusable launchers to demonstrate the novelty of the proposed procedure.**

## Nomenclature

| | | |
|---|---|---|
| $a$ | = | semi-major axis |
| $\mathbf{a}_D$ | = | aerodynamic drag acceleration vector |
| $b$ | = | cost coefficient |

[*] Assistant Professor, Department of Aerospace Engineering, (Former Postdoctoral Fellow in Daniel Guggenheim School of Aerospace Engineering at Georgia Institute of Technology).
[†] Dutton-Ducoffe Professor, Associate Professor, Daniel Guggenheim School of Aerospace Engineering, AIAA Senior Member; kokiho@gatech.edu (Corresponding Author)

$C$ = cost

$C_F$, $C_I$ = cost for rocket family and independent rocket

$\Delta v$ = velocity increment

$\Delta v_{ideal}$ = ideal velocity increment

$\Delta v_{loss}$ = velocity loss

$\Delta v_{r,ideal}$ = required ideal velocity increment

$e$ = eccentricity

$F$ = calibration factor for series production and operation

$f_{PS}$ = factor for cost savings

$f_u$ = cost factor for upgrade of reusable stage

$g_0$ = gravitational acceleration at sea level

$\mathbf{g}$ = gravity acceleration vector

$h$ = altitude

$h_{rb}$ = re-entry burn initiation altitude

$inc$ = inclination angle

$I_{sp}$ = specific impulse

$J$ = objective function

$j$ = $j^{th}$ launch

$k$ = cost coefficient

$L$ = learning factor

$m$ = mass

$m_0$, $m_L$ = gross lift-off mass and payload mass

$n_l$, $n_s$ = number of launches and number of stages

$n_r$ = number of certified reuses

$n_{rp}$ = number of productions of reusable stage

$\mathbf{r}$ = position vector

$t$ = time

$t_{co}$ = coasting time before boost-back burn

$t_{ff}$ = freefall time before re-entry burn

$T$ = thrust level

$TWR$ = thrust-to-weight ratio

$\mathbf{u}$ = thrust unit vector

$V_r$ = relative velocity

$V_{r,re}$ = relative velocity requirement after re-entry burn

$\mathbf{v}$ = velocity vector

$w$ = weighting factor

*Greeks*

$\alpha$ = angle of attack

$\varepsilon$ = structure ratio

$\eta$ = axial load

$\theta$ = pitch angle relative to the local horizontal

*Superscripts*

$*$ = optimal values

$\{p\}$ = values at $p^{th}$ phase

$(k)$ = values for $k^{th}$ launcher

*Subscripts*

$0$ = initial values

$c$ = values for compatible stage

$cc$ = values for comparable stage

$E$ = values in the ECEF frame

$e$ = values for engine

$f$ = final values

$fr$ = values for fairing

$i$ = values for $i^{th}$ stage ($i$ = 1, 2)

$l$ = values in the launch site

$lb$ = lower bound

$max$ = maximum values

$min$ = minimum values

$p$ = values for propellant

$r$ = values for the remaining stage

$s$ = values for structure

$t$ = target values

$ub$ = upper bound

## I. Introduction

OVER the past few decades, the development of launch vehicles has been generally led by the government with the primary concerns of mission success, national security, and technological advancements. However, with the advent of the new space era, numerous private companies have emerged as technologies have matured and the demand for space transportation has increased. It has naturally brought about a paradigm shift in space development from government-led to private company-led. Private companies strive to maximize profits by reducing the costs incurred to develop, manufacture, and operate launch vehicles. In consequence, saving costs has become crucial in the recent field.

Leveraging commonalities is one of the most profitable strategies when developing multiple rockets since the components, such as engines and stages, incur significant costs [1, 2]. That is, sharing compatible components enables manufacturers to operate the same product line for different rockets and avoid additional development for the components, resulting in cost and time savings. It also could be an efficient approach to the increasing demand for various space launch vehicles capable of carrying a wide variety of satellites, spacecraft, and cargo to space [3, 4]. In this context, product family design, whose commonality plays a key role, is a compelling approach for launcher manufacturers.

Product family design is a concept used to reduce costs and satisfy market niches in designing a group of compatible products whose commonality in functions and characteristics exists [5-7]. While it is widely implemented in diverse industries, it has not been considered seriously in rocket design due to the aforementioned characteristics of rocket development. Subsequent rockets in a conventional rocket family have been derived from the baseline model through minor or substantial modifications. This sequential approach is more efficient than an independent design. Yet, it does not fully take advantage of product family design and, thus, does not provide an optimal solution [8, 9], necessitating a product family-based simultaneous design process for multiple rockets. Notably, for the recent market trend in which multiple launch vehicles for different missions are concurrently developed, manufactured, and operated, the simultaneous design process could make significant synergies in terms of costs, time, and market needs.

Despite its effectiveness, only a few papers have explored the simultaneous design of multiple rockets because it has not been so long since the concurrent design of multiple rockets has become prominent. Delta IV and Atlas V families are the notable reference models that leverage the commonalities in a family for various missions, but the details of the design process are unreleased [10-12]. Hofstetter presents a framework analyzing commonality for portfolios of aerospace systems during an architecting phase and carries out a case study on the Saturn family [1]. Aliakbargolkar et al. propose a development strategy for a rocket family based on the methodology of architecting super-heavy launch vehicles from various points of view [2]. Although the technical details for a simultaneous design are not provided since these studies deal with a high-level decision process, these papers imply that a rocket family could benefit from sharing components across multiple rockets. Brevault et al. introduce multidisciplinary design optimization (MDO) approach for designing a rocket family composed of an RLV and ELV [13]. Even though their approach differs from this research, the idea, strategy, and approach for the design of a rocket family considering commonalities are notable. Insights can be found in the previous efforts toward simultaneous design methods for an aircraft family [9, 14, 15]. Given that commonalities are adequately involved in a design process, the approaches proposed in those papers are noteworthy.

The highly nonlinear nature of trajectory optimization makes it challenging to integrate multiple rockets into a single process without compromising reliability and accuracy. Accordingly, in this paper, we focus on only two disciplines - sizing and trajectory design - rather than considering the entire process in order to secure reliability and establish baselines in designing a rocket family simultaneously while abandoning excessive complexity. Motivated by the previous research, this paper proposes a sizing procedure for a two-member rocket family whose payload capabilities differ entirely. Instead of independent or sequential executions, the proposed procedure optimizes rockets in a single process to leverage commonalities fully. Assuming that the modifications are insignificant, the procedure is able to handle the stages in not only the same level but also the different levels between rockets (e.g., the first stage of Launcher 1 and the second stage of Launcher 2)

Sizing of a rocket refers to determining the stage weight of a multi-stage rocket that satisfies the mission requirements for given properties (e.g., propulsion, aerodynamic, and structure properties). The most critical factor in rocket sizing is estimating the accurate velocity losses. Therefore, the history of the research on rocket sizing, in short, lies in the advancement of the methods used to calculate velocity losses. While traditional approaches simplify the calculation with various assumptions, state-of-the-art research incorporates sizing and trajectory modules in a single loop to increase accuracy and optimality [16-22]. The proposed algorithms in the recent studies are reliable in that practical factors can be considered, yet they have a disadvantage in manipulation due to the highly nonlinear nature of trajectory optimization. For a single rocket, this is not a considerable concern; however, the issue makes it complex to exploit the existing algorithms in the rocket family sizing, indicating that a more tractable approach is required.

Despite the different goals, we can find an insightful idea for a tractable approach in [23]. The paper not only proposes a versatile approach for a preliminary design of a launcher but presents various perspectives on the applicable algorithms. Furthermore, the expenses of a rocket during its life cycle should be regarded as a top priority as the cost becomes more important to satisfy the increasing demands of stakeholders. A widely used objective function, maximizing the payload weight over the gross lift-off weight (GLOW) or minimizing the GLOW for a given payload weight, does not pursue the most cost-efficient design due to the modular characteristics of stages (the specific cost often varies by stage); alternatively, a cost-aware design approach is necessary.

In this work, an integrated sizing algorithm with a trajectory module for an individual rocket is presented and is expanded to simultaneous sizing for a two-member rocket family. The algorithm transforms sizing and trajectory optimizations into a single optimal control problem, thereby ensuring optimality and accuracy without iterative manipulations. Then, based on this approach, the commonality check, a process that determines whether commonalities between rockets exist, is performed by calculating the minimum and maximum weight of a compatible stage of each rocket; this step serves the role of determining whether the simultaneous sizing for rocket family is a feasible option and, if it is feasible, identifying that feasible design space so that the integrated optimizer can return the optimal solution within that space efficiently. Meanwhile, a simplified but realistic cost model for a rocket family is established with the consideration of commonalities. With the cost model, an optimal control problem for each rocket with the given compatible stage weight is formulated to optimize trajectory and the remaining stage weight that minimizes the specific cost (USD/kg). As a result, the simultaneous sizing for a rocket family becomes an optimization problem with single-dimensional search spaces over the compatible stage weight and is solved with a direct search method. A case study is carried out on a rocket family that involves an existing reusable launch vehicle (RLV), Falcon 9, and a preliminary designed expendable launch vehicle (ELV) to demonstrate the effectiveness of the proposed sizing procedure.

The contributions of this paper are as follows. First, an advanced sizing algorithm for a single rocket is proposed, which is more tractable due to non-iterative manipulation and improves accuracy and optimality in both sizing and trajectory perspectives. Moreover, the proposed algorithm is applicable to ELVs as well as RLVs. Second, a cost model that can capture the advantage of exploiting commonalities is established, enabling the costs to become a critical factor in rocket family design. Third, a proposed simultaneous sizing approach of a rocket family builds a baseline that possibly derives the succeeding research opportunities along with the growth of rocket industries. Furthermore, it broadens the design strategies manufacturers can take in planning multiple-launcher services.

## II. Unified sizing algorithm for single rockets

Sizing is a critical process in rocket design related directly to payload capability and costs. It determines the weight of each stage that meets the mission requirements for given properties, such as propulsion, structural, and aerodynamic features. Traditional sizing solves a simple optimization problem based on the Tsiolkovsky rocket equation with the estimations of velocity losses. Although this method is simple and thus easy to implement in the entire design process, the estimations entail inaccurate results since no analytic form for calculating velocity losses exists. Specifically, the inaccuracy becomes severe if the constraints imposed on the trajectory are intensive. To address this issue, recent studies utilize trajectory optimization in place of estimating the velocity losses and iteratively execute sizing and trajectory modules that interact with each other until convergence [18-21]. These methods are effective and reliable for a single rocket but still inappropriate to be applied to a simultaneous design of multiple rockets due to their iterative feature. As a solution to this matter, this section unifies and transforms the two modules into a single optimal control problem, eliminating the iterative executions and increasing accuracy and optimality.

### A. Review: Existing algorithms

Inputs: Mission requirements and aerodynamics, propulsion, and structure properties

**Optimization module**

$$\max J = m_L / m_0$$

Subject to

[Approach 1]

$$\sum_{i=1}^{n_s} \Delta v_i - \left( \Delta v_{r,ideal} + \Delta v_{loss} \right) = 0$$

$$=> m_{s_i} = m_{0_i} (1 - e^{-\frac{\Delta v_{tot_i}}{I_{sp\,i} g_0}}) \left( \frac{\varepsilon_i}{1-\varepsilon_i} \right)$$

[Approach 2]

$$\sum_{i=1}^{n_s} \Delta v_{ideal_i} - \Delta v_{r,ideal} = 0$$

$$=> m_{s_i} = m_{0_i} (1 - e^{-\frac{\Delta v_{ideal_i} + \Delta v_{loss_i}}{I_{sp\,i} g_0}}) \left( \frac{\varepsilon_i}{1-\varepsilon_i} \right)$$

**Velocity loss estimation**

[Method 1]
Historical data
[Method 2]
Simplified methods
[Method 3]
Trajectory optimization

Fig. 1 Summary of existing sizing methodologies

Fig. 1 summarizes the existing sizing methodologies, where $m_L$, $m_0$, $n_s$, $I_{sp}$, and $g_0$ are payload mass, GLOW, number of stages, specific impulse, and gravitational acceleration on the ground, respectively. Additionally, $\varepsilon$ is the structure ratio defined by structure (or dry) mass ($m_s$) and propellant mass ($m_p$) as $\varepsilon = m_s/(m_s + m_p)$, and the value is assumed to be given. The structure is assumed to include all components, subsystems, and propellant not consumed during the flight, such as residual and chill-down propellant. In practice, the structure ratio varies depending on the geometry of a launcher, such as diameter and length. Although this paper assumes that the structural properties are given, these values should be seriously analyzed in the entire design process since the impact of the structure ratio is considerable. $\Delta v$ and $\Delta v_{ideal}$ are, respectively, the velocity increment generated from a rocket propulsion system and the ideal (or actual) velocity increment achieved by a rocket. Moreover, $\Delta v_{r,ideal}$ and $\Delta v_{loss}$ are the velocity increment ideally required for payload insertion to the target orbit and velocity loss, respectively. The subscript $i$ indicates the values for the $i^{th}$ stage. Here, mission requirements (information on the target orbit and required payload mass, which is reinterpreted as the minimum payload mass), propulsion specifications (thrust and specific impulse), and aerodynamic and structure properties (aerodynamics coefficients, the diameter of a vehicle, nozzle exit area, and structure ratio) are assumed to be given.

Traditionally, the optimization approach with historical data or simplified methods for velocity loss estimations has been widely used. This approach distributes the required velocity increment ( $\Delta v_{r,ideal} + \Delta v_{loss}$ ) for the mission to each stage to maximize the specific payload weight ( $m_L/m_0$ ). $\Delta v_{r,ideal}$ is calculated from the orbital velocity of the target orbit and the Earth's rotation velocity on the launch site, and the velocity loss comes from historical data or simplified estimation methods. Regarding this approach, Jamilnia and Naghash adopted trajectory optimization for estimating the velocity loss to enhance accuracy and optimality [18]. However, despite the accurate velocity loss in each iteration, the first approach provides the results with a relatively large optimality gap since the total velocity increment is directly distributed to each stage without concerning the portions of the velocity loss in each stage.

An alternative approach developed by Koch [19] refines the above optimization approach by dividing it into two steps: 1) optimization without considering the velocity losses and 2) enlarging each stage by adding the velocity loss corresponding to the stage. What differentiates this approach from the first is that the optimization module in this approach determines only the ideal velocity increment for each stage. After finding the ideal velocity increment, the velocity loss for each stage obtained from trajectory optimization is added to calculate the stage mass, thereby serially enlarging the stages from top to bottom. In addition, other than the sizing algorithms for ELVs, Jo and Ahn proposed an optimal staging algorithm for an RLV based on Koch's method, and subsequently, they also presented a semi-integrated staging method with trajectory optimization that minimizes the life-cycle cost of an RLV [20, 22].

### B. Optimal control problem for sizing and trajectory optimization

The algorithms introduced in Subsection II-A cannot guarantee optimal solutions and are intractable in a simultaneous sizing process for multiple rockets. Accordingly, this subsection inserts the sizing module into the trajectory optimization by defining the mass values as a design variable and imposing additional constraints regarding stage properties. Consequently, the optimal control problem that can simultaneously handle sizing and trajectory is established as $P_0$. With the two optimal control problems whose set of design variables and constraints differs by the type of launcher, the proposed algorithm can cover the case of a family involving not only an ELV but also an RLV. For simplicity, this paper deals with only two-stage rockets.

**$P_0$: Optimal control problem for sizing and trajectory optimization of an ELV or RLV**

$$\min J = C / m_L \tag{1}$$

Subject to

Dynamic Constraints

$$\dot{\mathbf{r}} = \mathbf{v},\ \dot{\mathbf{v}} = \mathbf{g} + \mathbf{a}_D + \frac{T}{m}\mathbf{u},\ \dot{m} = -\frac{T}{I_{sp} g_0} \tag{2}$$

| Type of Rockets | ELV | RLV |
| --- | --- | --- |
| Design Variables | $\mathbf{u}(t),\ T^{\{1\}}(t),\ t_f^{\{p\}},\ m_{s_1},\ m_{s_2},\ m_L$ (3) | $\mathbf{u}(t),\ T^{\{1\}}(t),\ T^{\{6\}}(t),\ t_f^{\{p\}},\ m_{s_1},\ m_{s_2},\ m_L$ (9) |
| Boundary & Endpoints Constraints | $[\mathbf{r}_0^{\{1\}}, \mathbf{v}_0^{\{1\}}] = \left[\mathbf{r}_l, \mathbf{v}_l\right]$<br>$[a_f^{\{2\}}, e_f^{\{2\}}, inc_f^{\{2\}}] = \left[a_t, e_t, inc_t\right]$ (4) | $[\mathbf{r}_0^{\{1\}}, \mathbf{v}_0^{\{1\}}] = \left[\mathbf{r}_l, \mathbf{v}_l\right]$<br>$[a_f^{\{2\}}, e_f^{\{2\}}, inc_f^{\{2\}}] = \left[a_t, e_t, inc_t\right]$ (10)<br>$[\mathbf{r}_{E,f}^{\{6\}}, \mathbf{v}_{E,f}^{\{6\}}, \theta_f^{\{6\}}] = \left[\mathbf{r}_{E,t}, \mathbf{v}_{E,t}, 90^\circ\right]$<br>$m_f^{\{1\}} = m_0^{\{2\}} + m_0^{\{3\}}$ (11)<br>$[\mathbf{r}_0^{\{3\}},\ \mathbf{v}_0^{\{3\}}] = f_1\left(\mathbf{r}_f^{\{1\}},\ \mathbf{v}_f^{\{1\}}\right),\ t_0^{\{3\}} = t_f^{\{1\}} + t_{co}$<br>$[\mathbf{r}_0^{\{4\}},\ \mathbf{v}_0^{\{4\}}] = f_2\left(\mathbf{r}_f^{\{3\}},\ \mathbf{v}_f^{\{3\}}\right),\ t_0^{\{4\}} = t_f^{\{3\}} + t_{ff}$ (12)<br>$h_0^{\{4\}} = h_{rb},\ V_{r,f}^{\{4\}} \le V_{r,rb}$ |
| Path Constraints | $\left\lvert\alpha^{\{1\}}(t)\right\rvert \le \alpha_{a,\max},\ q^{\{1\}}(t) \le q_{a,\max}$ (5) | $\left\lvert\alpha^{\{1\}}(t)\right\rvert \le \alpha_{a,\max},\ q^{\{1\}}(t) \le q_{a,\max}$<br>$\left\lvert\alpha^{\{5\}}(t) - \pi\right\rvert \le \alpha_{d,\max},\ q^{\{5\}}(t) \le q_{d,\max}$ (13) |
| Constraints for Sizing | $TWR_0^{\{1\}} \le TWR_{\max},\ \eta_f^{\{1\}} \le \eta_{\max},\ \eta_f^{\{2\}} \le \eta_{\max}$ (6)<br>$\dfrac{m_{s_1}}{m_{s_1} + m_{p_1}} \equiv \dfrac{m_f^{\{1\}} - m_0^{\{2\}}}{m_0^{\{1\}} - m_0^{\{2\}}} = \varepsilon_1$ (7)<br>$m_L \equiv$<br>$m_f^{\{2\}} - \left(m_0^{\{2\}} - m_f^{\{2\}}\right)\dfrac{\varepsilon_2}{1-\varepsilon_2} - m_{fr} \ge m_{L,\min}$ (8) | $TWR_0^{\{1\}} \le TWR_{\max},\ \eta_f^{\{1\}} \le \eta_{\max},\ \eta_f^{\{2\}} \le \eta_{\max}$ (14)<br>$\dfrac{m_{s_1}}{m_{s_1} + m_{p_1}} \equiv \dfrac{m_f^{\{6\}}}{\dot{m}^{\{1\}} t_f^{\{1\}} + m_0^{\{3\}}} = \varepsilon_1$ (15)<br>$m_{p_2} \equiv m_0^{\{2\}} - m_f^{\{2\}} = \dot{m}^{\{2\}}\left(t_f^{\{2\}} - t_0^{\{2\}}\right)$ (16)<br>$m_L \equiv$<br>$m_f^{\{2\}} - \left(m_0^{\{2\}} - m_f^{\{2\}}\right)\dfrac{\varepsilon_2}{1-\varepsilon_2} - m_{fr} \ge m_{L,\min}$ (17) |

The objective function is to minimize the specific cost of the payload since the minimum GLOW does not always guarantee the minimum cost, especially in RLVs [9, 22]. Eq. (2) describes the point-mass dynamics in the Earth-Centered Inertial (ECI) frame, where $\mathbf{r}$, $\mathbf{v}$, $\mathbf{g}$, $\mathbf{a}_D$, and $\mathbf{u}$ are position, velocity, gravitational acceleration, drag acceleration, and thrust unit vectors, respectively, and $T$ is the magnitude of thrust. Superscript $\{p\}$ represents the values in the $p^{th}$ phase, and subscripts 0 and $f$ indicate the initial and final values, respectively. The design variables - thrust unit vector ($\mathbf{u}(t)$), thrust magnitude for the first phase ($T^{\{1\}}(t)$), final time of each phase ($t_f^{\{p\}}$), structure mass of each stage ($m_{s_i}$), and payload mass ($m_L$)- are determined by solving the problems. The thrust level in the 6$^{th}$ phase ($T^{\{6\}}(t)$) is an additional design variable for an RLV during its final landing phase, and the propellant mass ($m_{p_i}$) is calculated from the structure mass and structure ratio of each stage. It is important to note that payload mass is a design variable since the mission requirement can be satisfied by setting the target payload mass as the minimum value, and thus, the payload mass could become heavier.

For two-staged rockets, an RLV consists of 6 phases (ascent phases with the first and second stage burn/boost-back burn/re-entry burn/descent gravity turn/terminal landing burn), while the flight sequence of an ELV is divided into 2 phases (ascent phases with the first and second stage burn). Phase 1, the ascent phase of stage 1, includes a vertical lift-off for a certain amount of time followed by a gravity turn maneuver according to Eq. (5). This paper assumes that the stage separation and initiation of the second stage burn occur immediately after the end of the first stage burn, and the fairing separation is neglected.

The boundary, endpoint, and path constraints are imposed for trajectory optimization only. Eqs. (4) and (10) describe the boundary constraints for the states of an ELV and RLV, where subscripts $l$ and $t$ denote the values at the launch site, and target, respectively. As the landing site rotates with the Earth, the target states for landing are defined in the Earth-Centered Earth-Fixed (ECEF) frame with subscript $E$. Moreover, $a$, $e$, and $inc$ are semi-major axis, eccentricity, and inclination angle, respectively. In Eq. (10), $\theta$ is the pitch angle relative to the landing site, and it should be 90 degrees for a vertical landing. Eq. (11) links the mass of the separated stages between the ascent and descent phases, and Eq. (12) enforces the relationships that satisfy the mission requirements during the descent phase of an RLV, such as the coasting time ($t_{co}$) before the boost-back burn, freefall time ($t_{ff}$) followed by the re-entry burn, and initiation altitude ($h_{rb}$) and final velocity of the re-entry burn ($V_{r,rb}$). To avoid extensive transverse aerodynamic stress in the dense atmosphere, the permissible range of angle of attack (AOA) and the dynamic pressure is constrained, as shown in Eqs. (5) and (13). More details regarding the constraints and flight sequence on an RLV can be found in Ref. [22].

The thrust-to-weight ratio (*TWR*) and axial load ( $\eta$ ) need to be smaller than a specific value to avoid considerable atmospheric loads, and *TWR* should be larger than 1 to secure enough acceleration for lift-off under atmospheric pressure, as shown in Eqs. (6) and (14). As the mass of each stage becomes the design variable, which is set to be free, additional constraints are required to specify the relationships corresponding to the mass properties as described in Eqs. (7), (15), and (16). From the relationships, the structure ratio constraint is imposed in Eqs. (7) and (15). In addition, Eqs (8) and (17), where $m_{fr}$ is the fairing mass, define the payload mass requirement. Except for the constraints for sizing, it is optional to add or modify the constraints in accordance with the mission requirements or design principles.

By solving the optimal control problem, we can obtain both sizing and trajectory optimization results through a single execution regardless of the type of launcher, which is a promising feature of the proposed approach. However, this approach cannot be utilized directly in a simultaneous sizing procedure for a rocket family but is slightly modified. Furthermore, the independently optimized sizing results of each launcher obtained from this approach are compared with those from a simultaneous sizing procedure in the case study. Even though this approach by itself is unsuitable for the sizing of a rocket family, it is the most competent algorithm for the sizing of a single rocket.

The problems for rockets involving more than two stages or boosters can be easily defined by expanding the proposed formulations. As the number of stages increases, the additional phases should be considered with corresponding design variables (control inputs, time, and structure weight). The structure ratio constraint for the additional phase should also be imposed while the original constraints hold with only changing the phase numbers. Additionally, we can handle a launcher with boosters by separating the first phase into two phases: 1) the first phase until the booster separation, 2) the second phase from the booster separation to the end of the first stage burn. Assuming that the structure ratio and propulsion properties are given for the boosters, the size of the booster can be optimized by defining the booster weight as a design variable for the first phase. In this case, the structure ratio constraint for the booster is as follows

$$\frac{m_{s_b}}{m_{s_b}+m_{p_b}} \equiv \frac{m_f^{\{1\}}-m_0^{\{2\}}}{m_0^{\{1\}}-m_0^{\{2\}}-\dot{m}_{p_1}t_f^{(1)}}=\varepsilon_b, \tag{18}$$

where subscript $b$ denotes the values for the booster. The computational load will inevitably increase proportionally to the number of phases.

## III. Simultaneous sizing procedure for rocket families

This section contains the details of a simultaneous sizing procedure for rocket families proposed in this paper. The cost model for the rocket family is first defined to formulate the problem, followed by the design and optimization procedure.

### A. Cost Model

To formulate a simultaneous sizing problem for rocket families, a cost model needs to be developed that captures the cost savings by a common stage design shared among multiple rockets. Although the GLOW over the payload mass has been conventionally recognized as an objective function to be minimized to access a rocket design, minimizing GLOW does not necessarily lead to the lowest cost since the specific cost for each stage is not the same particularly when a rocket family is considered. In some cases, it could be more cost-efficient to enlarge a cheaper stage and downsize a more expensive stage, even with increasing GLOW. This mismatch between GLOW and cost is particularly severe in an RLV whose first stage is reused during its life cycle. Accordingly, this subsection introduces the cost model for a rocket family to find the cost-efficient mass configuration of rockets.

Based on the cost estimation relationships (CER), a parametric cost model is established to evaluate the life-cycle cost of a rocket family. Since it is derived from the existing models for an independent rocket developed in Refs.[22] and [24] and utilized in $P_0$ (only for the case study in this paper), the cost model for an independent rocket is introduced firstly as follows.

*Expendable launcher*

$$C_I = \sum_{i=1}^{n_s}\left(C_{D,fu_i} + C_{D,e_i}\right) + C_{D,fr} + \sum_{j=1}^{n_l}\left(\sum_{i=1}^{n_s}\left(C_{P,fu_i} + C_{P,e_i}\right) + C_{P,fr} + C_O\right)_j \tag{19}$$

*Reusable launcher*

$$\begin{aligned} C_I = &\sum_{i=1}^{n_s}\left(C_{D,fu_i} + C_{D,e_i}\right) + C_{D,fr} + \sum_{z}^{n_{rp}}\left(C_{P,fu_1} + C_{P,e_1}\right)_z \\ &+ \sum_{j=1}^{n_l}\left(\sum_{i=2}^{n_s}\left(C_{P,fu_i} + C_{P,e_i}\right) + C_{P,fr} + C_O\right)_j + \sum_{j=1}^{n_l-n_{rp}}\left(C_R\right)_j \end{aligned} \tag{20}$$

where, $C_I$ is the life-cycle cost of an independent rocket that consists of the development cost ($C_D$), production cost ($C_P$), operation cost ($C_O$), and cost for reuse ($C_R$), which is included only for an RLV, and the costs with subscripts *fu*, *e*, and *fr* correspond to fuselage, engine, and fairing, respectively. Note that the fuselage of a rocket represents the structure without an engine. Additionally, $n_s$ and $n_l$ are, respectively, the number of stages and launches. Moreover, $n_{rp}$ is the number of productions of the reusable stage determined as follows.

$$n_{rp} = \left\lceil \frac{n_l}{n_r + 1} \right\rceil, \tag{21}$$

where, $n_r$ is the certified (maximum) number of reuses. It is reasonable to assume that the reusable stage can be used a certain number of times. Hence, if the certified number of reuses is limited ( $n_r + 1 < n_l$ ), the reusable stage must be reproduced to achieve the target number of launches.

The model contains two formulations for both ELVs and RLVs, and each cost is categorized into development, production, and operation costs. As described in Eq. (19), in the case of an ELV, the development cost is nonrecurring, while the production and operation costs are recurring. On the other hand, the production cost of a reusable stage is nonrecurring, and the costs for reuse, such as refurbishment, transportation, and storage, are added in the RLV case. This model is expanded to the one for a rocket family by considering commonalities between the rockets to capture cost savings from sharing a compatible stage as follows. For simplicity, we assume two-stage rockets only in this paper.

*No commonality*

$$C_F = \sum_{k=1}^{2} C_I^{(k)} \tag{22}$$

*Commonality existence*

<u>*1) Expendable + Expendable*</u> ( $m_c$ : $m_{s_1}^{(1)} = m_{s_1}^{(2)}$ , $m_{s_2}^{(1)} = m_{s_2}^{(2)}$ , *or* $m_{s_2}^{(1)} = m_{s_1}^{(2)}$ )

$$\begin{aligned} C_F^{(1)} &= C_{D,fu_r}^{(1)} + \sum_{i=1}^{2}\left(C_{D,e_i}^{(1)}\right) + C_{D,fr}^{(1)} + \sum_{j=1}^{n_l^{(1)}}\left(C_{P,fu_r}^{(1)} + \sum_{i=1}^{2}\left(C_{P,e_i}^{(1)}\right) + C_{P,fr}^{(1)} + C_O^{(1)}\right)_j \\ &\quad + \left(C_{D,fu_c} + f_{PS}\sum_{j=1}^{n_l^{(1)}+n_l^{(2)}}\left(C_{P,fu_c}\right)_j\right)\frac{n_l^{(1)}}{n_l^{(1)}+n_l^{(2)}}, \end{aligned} \tag{23}$$

$$\begin{aligned} C_F^{(2)} &= C_{D,fu_r}^{(2)} + \sum_{i=1}^{2}\left(C_{D,e_i}^{(2)}\right) + C_{D,fr}^{(2)} + \sum_{j=1}^{n_l^{(2)}}\left(C_{P,fu_r}^{(2)} + \sum_{i=1}^{2}\left(C_{P,e_i}^{(2)}\right) + C_{P,fr}^{(2)} + C_O^{(2)}\right)_j \\ &\quad + \left(C_{D,fu_c} + f_{PS}\sum_{j=1}^{n_l^{(1)}+n_l^{(2)}}\left(C_{P,fu_c}\right)_j\right)\frac{n_l^{(2)}}{n_l^{(1)}+n_l^{(2)}} \end{aligned} \tag{24}$$

*2) Reusable + Expendable* ( $m_c$ : $m_{s_2}^{(1)} = m_{s_2}^{(2)}$ or $m_{s_2}^{(1)} = m_{s_1}^{(2)}$ )

$$\begin{aligned} C_F^{(1)} = {} & C_{D,fu_1}^{(1)} + \sum_{i=1}^{2}\left(C_{D,e_i}^{(1)}\right) + C_{D,fr}^{(1)} + \sum_{z=1}^{n_{rp}^{(1)}}\left(C_{P,fu_1}^{(1)} + C_{P,e_1}^{(1)}\right)_z + \sum_{j=1}^{n_l^{(1)}}\left(C_{P,e_c}^{(1)} + C_{P,fr}^{(1)} + C_O^{(1)}\right)_j + \sum_{j=1}^{n_l^{(1)}-n_{rp}^{(1)}}\left(C_R\right)_j \\ & + \left(C_{D,fu_c} + f_{PS}\sum_{j=1}^{n_l^{(1)}+n_l^{(2)}}\left(C_{P,fu_c}\right)_j\right)\frac{n_l^{(1)}}{n_l^{(1)}+n_l^{(2)}}, \end{aligned} \tag{25}$$

$$\begin{aligned} C_F^{(2)} = {} & C_{D,fu_r}^{(2)} + \sum_{i=1}^{2}\left(C_{D,e_i}^{(2)}\right) + C_{D,fr}^{(2)} + \sum_{j=1}^{n_l^{(2)}}\left(C_{P,fu_r}^{(2)} + \sum_{i=1}^{2}\left(C_{P,e_i}^{(2)}\right) + C_{P,fr}^{(2)} + C_O^{(2)}\right)_j \\ & + \left(C_{D,fu_c} + f_{PS}\sum_{j=1}^{n_l^{(1)}+n_l^{(2)}}\left(C_{P,fu_c}\right)_j\right)\frac{n_l^{(2)}}{n_l^{(1)}+n_l^{(2)}} \end{aligned} \tag{26}$$

*3) Reusable + Reusable* ( $m_c$ : $m_{s_2}^{(1)} = m_{s_2}^{(2)}$ )

$$\begin{aligned} C_F^{(1)} = {} & C_{D,fu_1}^{(1)} + \sum_{i=1}^{2}\left(C_{D,e_i}^{(1)}\right) + C_{D,fr}^{(1)} + \sum_{z=1}^{n_{rp}^{(1)}}\left(C_{P,fu_1}^{(1)} + C_{P,e_1}^{(1)}\right)_z + \sum_{j=1}^{n_l^{(1)}}\left(C_{P,e_c}^{(1)} + C_{P,fr}^{(1)} + C_O^{(1)}\right)_j + \sum_{j=1}^{n_l^{(1)}-n_{rp}^{(1)}}\left(C_R^{(1)}\right)_j \\ & + \left(C_{D,fu_c} + f_{PS}\sum_{j=1}^{n_l^{(1)}+n_l^{(2)}}\left(C_{P,fu_c}\right)_j\right)\frac{n_l^{(1)}}{n_l^{(1)}+n_l^{(2)}}, \end{aligned} \tag{27}$$

$$\begin{aligned} C_F^{(2)} = {} & C_{D,fu_1}^{(2)} + \sum_{i=1}^{2}\left(C_{D,e_i}^{(2)}\right) + C_{D,fr}^{(2)} + \sum_{z=1}^{n_{rp}^{(2)}}\left(C_{P,fu_1}^{(2)} + C_{P,e_1}^{(2)}\right)_z + \sum_{j=1}^{n_l^{(2)}}\left(C_{P,e_c}^{(2)} + C_{P,fr}^{(2)} + C_O^{(2)}\right)_j + \sum_{j=1}^{n_l^{(2)}-n_{rp}^{(2)}}\left(C_R^{(2)}\right)_j \\ & + \left(C_{D,fu_c} + f_{PS}\sum_{j=1}^{n_l^{(1)}+n_l^{(2)}}\left(C_{P,fu_c}\right)_j\right)\frac{n_l^{(2)}}{n_l^{(1)}+n_l^{(2)}} \end{aligned} \tag{28}$$

where, $C_F$ is the life cycle cost of a rocket family, and subscripts *c* and *r* indicate the compatible and remaining stages, respectively. Note that the other stage is defined as a remaining stage once the compatible stage is specified in each rocket. Additionally, $f_{PS}$ is a production-sharing factor that reflects the cost savings resulting from sharing the same product line.

If commonality does not exist, the cost for a family is simply calculated by adding up the independent costs. On the contrary, when commonality is found, the cost is evaluated with respect to the three cases: 1) all launchers in a family are expendable, 2) one is reusable, and the other is expendable, and 3) all launchers are reusable. A reusable stage is not regarded as a compatible stage in this paper. Of course, technically, a reusable stage could be shared, but in this case, the effect of cost savings may be limited (or the cost becomes more expensive) because the nonrecurring cost of manufacturing a reusable stage leads to the trivial benefits of sharing the product line of a reusable stage.

As the two launchers utilize the same stage, a single development of the compatible stage can cover the two launchers, eliminating the development cost and period for the stage of the second launcher. Furthermore, as the number of manufacturing for the stage increases, the production cost takes full advantage of the learning effect and production facility sharing. Accordingly, the development and production costs for a compatible part are calculated from a single equation and distributed to each launcher in proportion to the number of launches.

We can express the cost equations in Eqs. (23)-(28) in a different approach by defining the cost savings ( $C_{savings}$ ) as follows

$$C_F^{(k)} = C_I^{(k)} - C_{savings}^{(k)}, \tag{29}$$

$$C_{savings}^{(k)} = C_{D,fu_c} + \sum_{j=1}^{n_l^{(1)}} \left( C_{P,fu_c} \right)_j - \left( C_{D,fu_c} + f_{PS} \sum_{j=1}^{n_l^{(1)}+n_l^{(2)}} \left( C_{P,fu_c} \right)_j \right) \frac{n_l^{(1)}}{n_l^{(1)} + n_l^{(2)}}. \tag{30}$$

We can calculate the saved costs thanks to sharing the compatible part from Eq. (30).

Assuming that the cost is a function of mass only, the cost elements are defined as [22, 24]

$$C_D = k_D m_{ref}^{b_D},\ C_P = k_P m_{ref}^{b_P} F_P,\ C_O = k_O m_{ref}^{b_O} F_O,\ C_R = k_R m_{ref}^{b_R} F_R, \tag{31}$$

where, $m_{ref}$ is the reference mass corresponding to each component, and $k$ and $b$ are the coefficients defined through the cost analysis on a component versus mass. $F$ is the calibration factor for series production and operation described as

$$F = j^{\log_2 L}, \quad (j = 1, \ldots, n_l). \tag{32}$$

Through Eq. (32) with the learning factor $L$, the learning effect is applied when serially performing the same tasks. Generally, the production and operation costs for the $j^{th}$ launch decrease compared to the previous one due to the learning effect, where the manufacturer becomes more efficient in handling tasks. In contrast, the cost for reuse, a type of operation cost, may become more expensive to deal with the fatigue of a reusable stage as it is repeatedly used. In practice, this cost increment should be initialized when we start producing a new unit for the reusable stage, but this paper does not consider this fact in the cost model for simplicity and alternatively selects the appropriate value for the learning factor.

### B. Design procedure

With the cost function defined, the rocket family sizing problem can be formulated as a simultaneous optimization of multiple rockets with common variables; this problem is challenging because (1) there is no analytical derivative information due to the embedded optimal control problem; (2) the constraints are not analytically represented and their evaluation involves solving optimal control problems; and (3) the feasibility space is often relatively small. Although this first challenge can be resolved using a finite difference approach or a derivative-free direct search algorithm such as the Nelder-Mead method [25, 26], it is challenging (or sometimes impossible) to ensure this type of solvers to find a feasible solution for highly-constrained problems. Therefore, we need an alternative approach that can ensure the feasibility of the solution. Furthermore, in practice, we often encounter real-world situations where we do not even know if the rocket family design problem itself has a feasible solution or not, so there is a need for checking the feasibility of the problem before running the full optimization.

In response to the identified challenges, we propose an alternative approach: first identifying the feasible space in the form of the bounds of a key parameter (i.e., commonality check) and then finding the solution within the bounds via direct search. First, the commonality check is performed to determine whether a common stage exists between the two rockets. It solely compares the weight between the two stages of different rockets (e.g., the first stage of launcher 1 and the second stage of launcher 2) by assuming that modifications are negligible when sharing the common stage. Concluded that the two stages are compatible, an optimization problem that minimizes the cost of a family is solved via a direct search method. The optimization problem consists of two optimal control problems as subproblems that find the remaining stage mass, payload mass, and trajectory of each launcher for the given compatible stage weight with the objective function of minimizing the specific cost of the payload.

Note that even though the commonality check that provides the bounds seems redundant from the perspective of computational efficiency in that it solves four additional optimal control problems, this step is critical because (1) it ensures that the final solution is feasible; and equally importantly, (2) it determines whether simultaneous optimization is necessary in the first place.

Fig. 2 illustrates a conceptual diagram of the proposed procedure, each process in this figure will be explained in more detail in the following subsubsections.

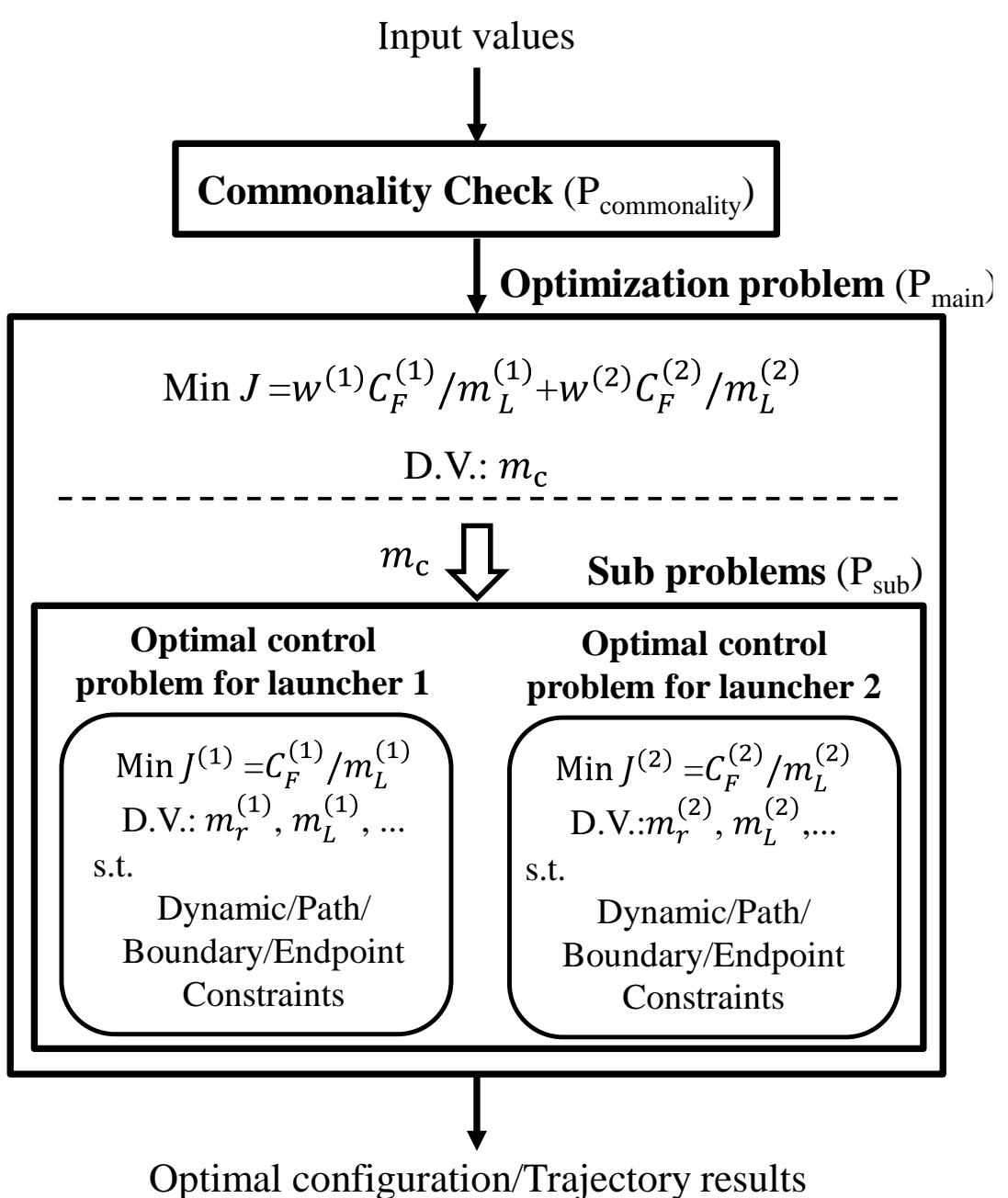

Fig. 2. Conceptual diagram of simultaneous sizing procedure for a rocket family

*1. Commonality Check*

The goal of the commonality check is to determine if commonality exists in a rocket family by finding the minimum and maximum weight of a comparable stage of each stage. At the beginning of a design process, the existence of comparable stages between two rockets is roughly examined through empirical data or simple calculations. In this paper, the comparable stage is defined as a candidate possibly compatible with the one in another rocket. Therefore, four cases of stage sets appear in a two-member family, each with two stages, and technically, the commonality check needs to be performed on each case to define the compatibility of a specific stage between the rockets. However, it is possible to empirically distinguish the case with the possibility of commonality when the required payload mass and properties of the rocket are given for each mission. If concluding that a specific stage in Launcher 1 is qualitatively comparable to the one in Launcher 2, the quantitative commonality check is performed with those stages. Fig. 3 illustrates the example of the commonality check: two launchers are considered for entirely different two missions, and it is evident that the potential commonality can only be observed between the second stage of Launcher 1 and the first stage of Launcher 2.

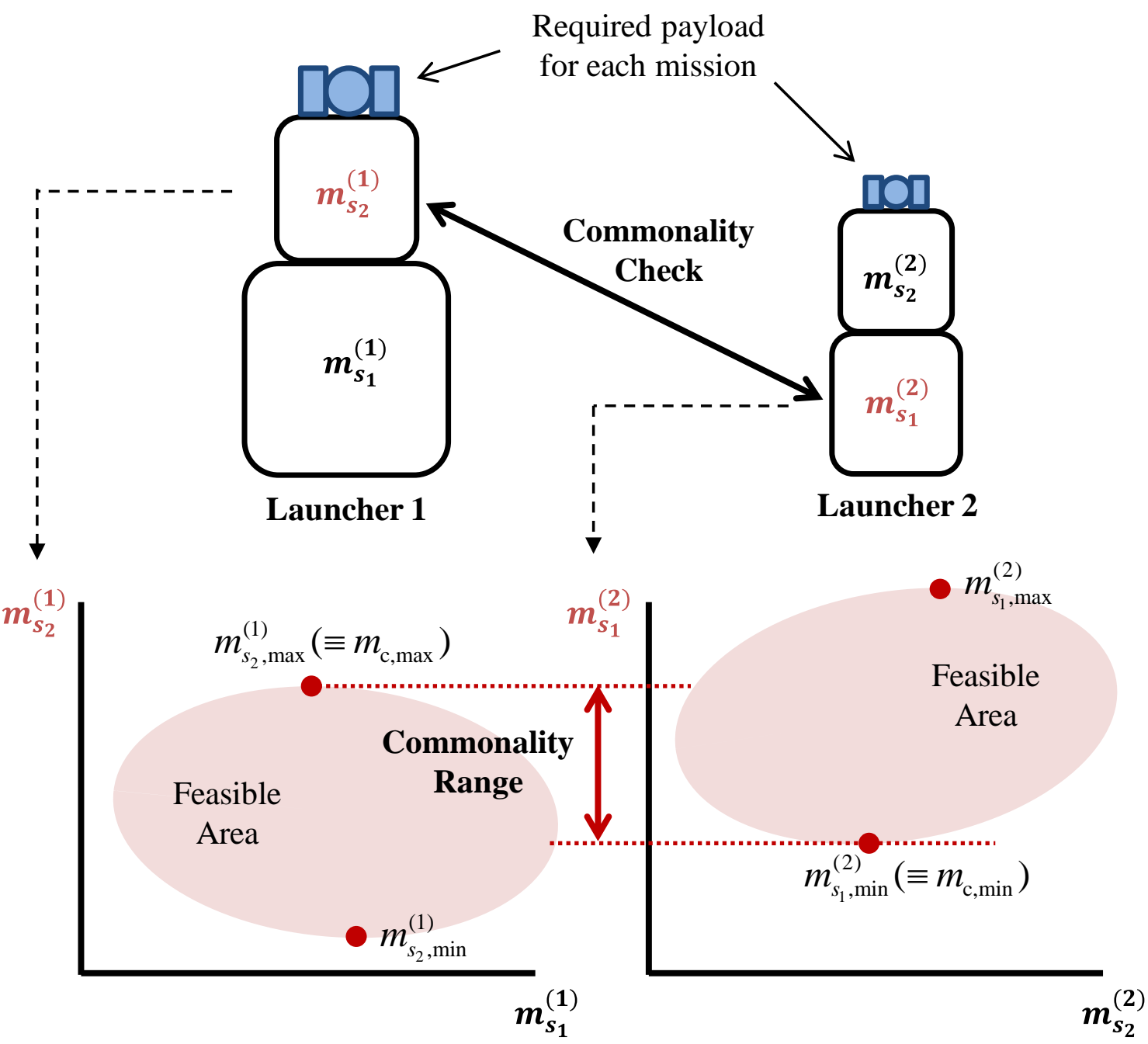

Fig. 3 Example of commonality check

The feasible area is a set of mass configurations capable of carrying a payload heavier than the required weight to the target orbit. Therefore, if the feasible area of one launcher overlaps with the other's, as shown in Fig 3, the possibility of commonality is inherent in the design process. However, the feasible area does not necessarily need to be obtained. Instead, after selecting the two stages comparable in terms of weight from the empirical data, each stage's maximum and minimum mass in the feasible area is used to evaluate the commonality between the two stages. That is, the feasible range over the comparable stage is calculated from its min/max mass, and if there is an overlap between the ranges, we can conclude that the two stages are compatible. Therefore, this step finds the minimum and maximum mass of the comparable stage for each launcher ($[m_{cc,\min}^{(1)}, m_{cc,\max}^{(1)}, m_{cc,\min}^{(2)}, m_{cc,\max}^{(2)}]$) through the slightly modified optimal control problems from $P_0$, defined as $P_{commonality}$. It is noteworthy that comparable stage ($m_{cc}$) indicates the stage with the possibility of being shared, whose commonality check is required and becomes the compatible stage ($m_c$) once the commonality is confirmed.

**$P_{commonality}$: Optimal control problems for min/max mass of the comparable stages (4 independent problems)**

| | Launcher 1 | | Launcher 2 | |
|---|---|---|---|---|
| Objective functions | $\min J = m_{cc}^{(1)}$ | $\max J = m_{cc}^{(1)}$ | $\min J = m_{cc}^{(2)}$ | $\max J = m_{cc}^{(2)}$ |
| Dynamic Constraints | Eq. (2) | | | |
| Design Variables | ELV: Eq. (3)<br>RLV: Eq. (9) | | | |
| Constraints | ELV: Eqs. (4)-(8)<br>RLV: Eqs. (10)-(17) | | | |

The objective function of each optimal control problem is to minimize or maximize the comparable stage mass subject to constraints identical to those described in $P_0$. Depending on the type of rocket, the corresponding design variables and constraints are selected for the problems: if the rocket is expendable, Eqs. (3)-(8) are used or Eqs. (9)-(17) are used for RLVs. Consequently, by supposing the launcher with a larger value for the maximum comparable stage mass is numbered launcher 1, the commonality range is determined by $m_{c,\max}$ and $m_{c,\min}$ which are calculated from the simple logic as follows.

$$\begin{cases} m_{c,\max} = m_{cc,\max}^{(2)},\ m_{c,\min} = m_{cc,\min}^{(1)}, & \text{if } m_{cc,\max}^{(1)} \ge m_{cc,\max}^{(2)} \ge m_{cc,\min}^{(1)}\ \&\ m_{cc,\min}^{(1)} \ge m_{cc,\min}^{(2)}, \\ m_{c,\max} = m_{cc,\max}^{(2)},\ m_{c,\min} = m_{cc,\min}^{(2)}, & \text{if } m_{cc,\max}^{(1)} \ge m_{cc,\max}^{(2)}\ \&\ m_{cc,\min}^{(1)} \le m_{cc,\min}^{(2)}, \\ \textit{No commonality}, & \textit{else}. \end{cases} \tag{33}$$

*2. Optimization for simultaneous sizing*

With the identified commonality bounds, this subsection provides the optimization approach in the simultaneous sizing of a rocket family. The considered optimization problem consists of two optimal control problems to be solved for the given compatible stage mass at every iteration loop and finds the optimal compatible stage mass. To this end, the compatible stage mass is set as a design variable with the objective function of minimizing the summation of the specific cost of each payload. In this regard, a weighting factor ($w$) is used to evaluate the objective function as the following equation because the payload mass is valued differently in each mission.

$$\min J = w^{(1)} C_F^{(1)} / m_L^{(1)} + w^{(2)} C_F^{(2)} / m_L^{(2)}. \tag{34}$$

The weighting factor should be carefully determined based on the importance of the missions, the type of target orbit, and the number of launches. For example, if the target orbit of Launcher 1 is a geosynchronous equatorial orbit (GEO) and the other one is a low-Earth orbit (LEO), then the weight factor of Launcher 1 should be considered more significantly, or the specific cost for the launcher, expected to have more launches, needs to be given a higher value. Despite the weighting factors, we selected this objective function since it can adequately capture the benefits of the increments in the payload weight. Leveraging a particular configuration that can increase the payload capacity at a low-cost increment is a profitable strategy, fulfilling a more comprehensive range of customer needs. If the cost is of utmost importance, we can select the objective function minimizing the total cost ( $\min J = C_F^{(1)} + C_F^{(2)}$ ) instead of the current form. The main optimization problem is formulated as follows:

**$P_{main}$: Main optimization problem**

| | Simultaneous optimization |
|---|---|
| Objective functions | $\min J = w^{(1)} C_F^{(1)} / m_L^{(1)} + w^{(2)} C_F^{(2)} / m_L^{(2)}$ |
| Design Variables | $m_c$ |
| Subproblems | $P_{sub} \Rightarrow \begin{matrix} C_F^{(1)} / m_L^{(1)} \ \textit{for Launcher}\ 1 \\ C_F^{(2)} / m_L^{(2)} \ \textit{for Launcher}\ 2 \end{matrix}$ |

The optimal control problems, the subproblems of the optimization problem for simultaneous sizing, are defined as follows.

**$P_{sub}$: Optimal control subproblems for given compatible stage mass (2 independent problems)**

| | Launcher 1 | Launcher 2 |
|---|---|---|
| Objective functions | $\min J = C_F^{(1)} / m_L^{(1)}$ | $\min J = C_F^{(2)} / m_L^{(2)}$ |
| Dynamic Constraints | Eq. (2) | |
| Design Variables | Case 1. ELV ( $m_c = m_{s_1}$ ): | $\mathbf{u}(t),\ t_f^{\{2\}},\ m_{s_2} (\equiv m_r),\ m_L$ (35) |
| | Case 2. ELV ( $m_c = m_{s_2}$ ): | $\mathbf{u}(t),\ t_f^{\{p\}},\ m_{s_1} (\equiv m_r),\ m_L$ (36) |
| | Case 3. RLV ( $m_c = m_{s_2}$ ): | $\mathbf{u}(t),\ T^{\{6\}}(t),\ t_f^{\{p\}},\ m_{s_1} (\equiv m_r),\ m_L$ (37) |
| Constraints | Case 1. ELV ( $m_c = m_{s_1}$ ): | Eqs. (4)-(6) and (8)<br>$m_f^{(1)} - m_0^{(2)} = m_c$ (38) |
| | Case 2. ELV ( $m_c = m_{s_2}$ ): | Eqs. (4)-(7)<br>$m_L \equiv m_f^{\{2\}} - m_c - m_{fr} \geq m_{L,\min}$ (39) |
| | Case 3. RLV ( $m_c = m_{s_2}$ ): | Eqs. (10)-(15)<br>$m_L \equiv m_f^{\{2\}} - m_c - m_{fr} \geq m_{L,\min}$ (40) |

$P_{sub}$ encompasses the two independent optimal control problems whose objective function is to minimize the specific cost of the payload subject to the design variables and constraints defined by which stage is compatible and the type of launcher. Any stage of an ELV can be compatible, while only the second stage of an RLV can, and thus three cases exist with distinct problem formulations for each case. The same dynamic constraints are used, and several design variables and constraints are identical to those in $P_0$, but few are slightly modified as the compatible stage mass ( $m_c$ ), a design variable in $P_0$, is given. As a result, only the mass of the remaining stage ( $m_r$ ) between the two stages remains as a design variable. Moreover, if the first stage in an ELV is a compatible stage, which is given, then the final time of phase 1 is automatically determined, and the additional constraint in Eq. (38) is imposed to handle the discontinuity between the two phases caused by the stage separation. In contrast, in Eqs. (36) and (37), the time values for phases 1 and 2 are design variables despite the given second stage mass because the time for the first phase is unknown. Additionally, a minor adjustment is made in calculating the payload mass in Eq. (17) as described in Eqs. (39) and (40).

In $P_{main}$, the performance index to be minimized is a function of $m_c$ as $f(m_c) = w^{(1)} C_F^{(1)} / m_L^{(1)} + w^{(2)} C_F^{(2)} / m_L^{(2)}$ in which $f(m_c)$ is evaluated by solving the optimal control problems $P_{sub}$. Consequently, the sizing of a rocket family is transformed into an optimization problem $P_{main}$ with a single-dimensional search space over the compatible stage mass $m_c$ , whose feasible range is bracketed by the upper and lower bounds found in the commonality check step. As mentioned previously, a direct search method, a technique for solving the optimization problem without exploiting the gradient of the objective function, is employed to solve the problem $P_{main}$. Applying two gradient-based algorithms (sqp and interior point) with the initial guess within the commonality range, we observed that these algorithms are sensitive to the initial guess and occasionally fail to converge to a solution. A direct search method determines the best points for the next search by sequentially examining the trial solutions corresponding to the current search points and comparing them. Numerous algorithms have been developed so far [27], and among them, this paper presumes the golden section search is most suitable for this problem in that it is straightforward and easy to implement, the bounds for search spaces are provided from the commonality check, and the problem has only one design variable [25]. Notably, it is imperative to search for a solution within the bounds in which the compatible mass guarantees the feasibility for both launchers.

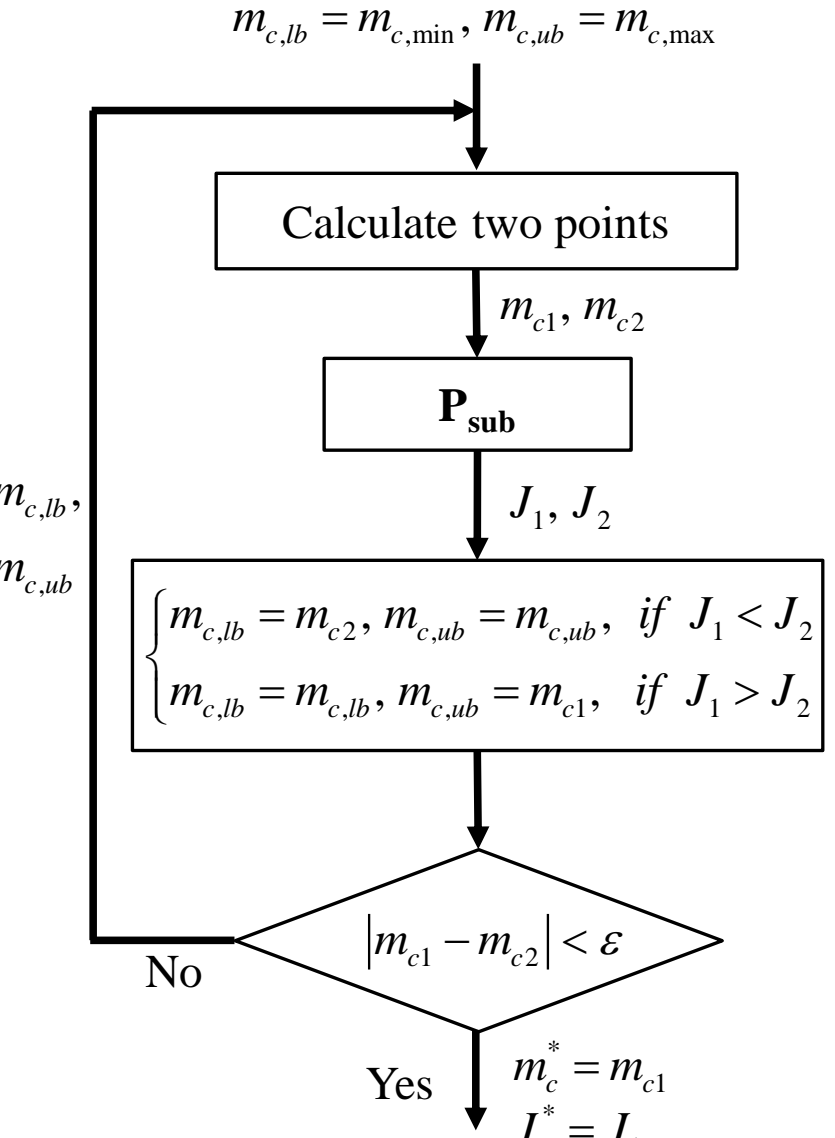

Fig. 4. Flow chart of the optimization problem for simultaneous sizing with the golden section search method

Fig. 4 depicts the flow chart of the optimization problem in simultaneous sizing for a rocket family with the golden section search method. The method starts with calculating the two search points with the golden ratio ($\varphi$) of 1.618 as follows [25].

$$\begin{aligned} m_{c1} &= m_{c,lb} + (\varphi - 1)(m_{c,ub} - m_{c,lb}), \\ m_{c2} &= m_{c,ub} - (\varphi - 1)(m_{c,ub} - m_{c,lb}), \end{aligned} \tag{41}$$

where, $m_{c,lb}$ and $m_{c,ub}$ are the lower and upper bounds of the compatible stage mass defined as $m_{c,ub} = m_{c,\max}$ and $m_{c,lb} = m_{c,\min}$, respectively. Afterward, the objective function value in Eq. (34) is evaluated by solving $P_{sub}$ for each search point. With the comparison of the values, the bounds are then updated based on the following logic.

$$\begin{cases} m_{c,lb} = m_{c2},\ m_{c,ub} = m_{c,ub}, & \text{if } J_1 < J_2, \\ m_{c,lb} = m_{c,lb},\ m_{c,ub} = m_{c1}, & \text{if } J_1 > J_2. \end{cases} \tag{42}$$

As this logic aims to find the minimum value, if the value for a specific search point is confirmed to be greater than the other between the two, there is no need to continue searching beyond that point. Therefore, this particular point, associated with the greater function value, becomes the bound according to Eq. (42). This process continues until the two search points converge. Then, as the solutions of the procedures, the compatible stage mass is directly optimized through the golden section method, and the remaining stage mass with trajectory results is obtained from the optimal control problems in the last iteration loop.

Lastly, we would like to note that the golden section method is only guaranteed to find a global optimal solution for a function that does not have a local optimum, such as a convex or monotonic function. Although we cannot rigorously prove that convexity property for our problem due to the complexity behind the optimal control solvers involved in the trajectory design subproblems, we can still reasonably assume it for the following reasons. Within the feasible range of $m_c$, we have multiple tradeoffs occurring in the trajectory simulation: namely, a larger $m_c$ means a greater compatible stage performance, which can indicate a smaller remaining stage size ( $m_r$ ) to deliver the same payload; this trend holds until the stage mass $m_c$ becomes too large and the gravity-driven velocity loss (and steering loss if the $m_r$ is the second stage) becomes so dominant that larger $m_c$ starts to degrade the rocket performance, in which case a heavier remaining stage is needed to compensate for this performance loss.[‡] Considering these tradeoffs, $m_r$ is expected to follow a convex (or monotonic) function of $m_c$ within its feasible range. Although the actual cost function is nonlinear, as we observe a nearly linear trend within the realistic rocket mass range, we can reasonably assume that the total cost (summation of the cost for $m_r$ and the cost for $m_c$ ) is also convex (or monotonic) with respect to $m_c$. This above discussion is based on a qualitative argument, but this assumption is demonstrated to be valid even with realistic simulations and constraints as shown in Sec. IV, where the golden section method is able to find the global optimum for the tested cases.

Although this paper focuses on a two-member family for simplicity, it might be beneficial to discuss the potential application of the proposed methodology to a rocket family with more than two members. If all launchers have just a single compatible stage, the size of the problem (i.e., the number of optimal control problems to be solved) becomes linearly larger proportional to the number of members. In contrast, a launcher with more than one compatible stage increases the dimension of the problem size. For instance, in a three-member family whose the first and second stages of Launcher 2 are respectively compatible with Launcher 1 and 3, the size of $P_{commonality}$ doubles, and the size of $P_{main}$ becomes squared. This increment affects the computational load in a similar manner.

[‡] Note that the thrust and $I_{sp}$ are kept constant in this comparison, so a larger stage mass corresponds to a longer burn time and thus, a greater gravity loss.

## IV. Case Study

This section presents the case study results on a rocket family comprising two launch vehicles to prove the effectiveness of the proposed procedure. We select an existing RLV, Falcon 9, and a preliminary designed ELV aligned with the RLV as the family members. The details on the specifications and mission information are provided in Table 1 [28].

Table 1: Specifications of a launcher family

| | | Launcher 1 (RLV) | Launcher 2 (ELV) |
|---|---|---|---|
| Thrust, tonf | | 838. 8 (768.6), 95.2 | 240 (228.3), 40 |
| Specific impulse, s | | 311, 348 | 315, 340 |
| Structure ratio, % | | 5.13, 3.59 | 3.59, 4.5 |
| Fairing mass, kg | | 1,700 | 500 |
| Minimum payload mass, kg | | 15,628 | 5,766 |
| Target *a, e, inc,* km, -, deg | | 6,635.7, 0.0073, 51.6 | 6,678.2, 0.0001, 51.6 |
| Total number of launches | Case 1 | 50 ($n_r$ = 15) | 50 |
| | Case 2 | 150 ($n_r$ = 15) | 50 |
| | Case 3 | 50 ($n_r$ = 15) | 150 |
| | Case 4-A | 150 ($n_r$ = 50, $f_u$=1) | 50 |
| | Case 4-B | 150 ($n_r$ = 50, $f_u$=1.2) | 50 |

The values in each cell for the specification correspond to the first and second stages, respectively, and the value in parentheses represents the thrust at sea level. The information regarding aerodynamics, which is publicly unreleased, is estimated based on available empirical data. Each launcher in the family should be able to carry a payload of more than 15,628 kg and 5,766 kg to two different low Earth orbits in each mission, respectively. Although the target orbits are similar, the payload requirement differs entirely between the launchers, necessitating the two types of launchers to achieve the missions. Thus, we assume that the two launchers are designed to have comparable stages with an identical structure ratio, i.e., the second stage of Launcher 1 and the first stage of Launcher 2.

Table 2 represents the parameters used for the constraints in the optimal control problems. For simplicity, we assume that the thrust level for all phases is fixed but take into account pressure loss caused by atmospheric pressure. The absolute values for the maximum AOA during phases 1 and 5 are limited to 1 and 2 degrees, respectively. These values might be large for an actual flight but not lead to a significant variation from the trajectory results with a smaller value close to zero where convergence issues may occur. In this paper, we manage the maximum dynamic pressure through the TWR constraint and the predesigned values, $h_{rb}$ and $V_{r,rb}$ instead of directly imposing the max Q constraints. Since the max Q is proportional to TWR in principle, the maximum TWR constraint effectively controls the dynamic pressure constraint. Therefore, the maximum TWR is determined by analyzing dynamic pressure with various trajectory simulations to ensure the rockets ascend under an acceptable load. In determining the TWR value, we set the max Q as 40kPa, which is somewhat larger than the conventional value by assuming that it can be managed by throttling in an actual flight, even though throttling is not directly considered in this trajectory optimization.

Although Falcon 9 is equipped with grid fins to control/leverage the aerodynamic forces, especially to generate the lift force for the effective flight to the downrange direction, this case study does not consider the grid fins and assumes that the lift force is negligible. Technically, it is more reasonable to select the smaller value for the maximum axial load than 6g, yet this paper aims to demonstrate the broad area of the feasible region with a sufficiently large value, even though it is an inactive constraint in the optimal solutions. The parameters in Eq. (12) come from Ref. [22], where further details are provided. We exploited the golden-section search algorithm for $P_{main}$ as described in Subsection III-B and solved the optimal control problems in $P_{commonality}$ and $P_{sub}$ through GPOPS-II in MATLAB environments [29].

Table 2. Parameters for constraints

| Parameter | Value |
|---|---|
| $\alpha_{a,\max}$ , $\alpha_{d,\max}$ , deg | 1, 2 |
| $TWR_{\max}$ | 1.45 |
| $\eta_{\max}$ , g | 6 |
| $t_{co}$ , s | 20 |
| $h_{rb}$ , m | 800 |
| $V_{r,rb}$ , m/s | 60 |

## A. Case study 1: single rocket

This subsection presents case study results on a single launch vehicle to validate the effectiveness of the unified sizing algorithm proposed in Section II since it is one of the core contributions of this research. We applied the algorithm to Launcher 2 in Table 1 and compared the results with those obtained by utilizing the state-of-the-art method in [19]. In this case study, the objective function is to minimize the GLOW over the payload weight ($m_0/m_L$), which is the traditional performance index used for rocket sizing. While imposing constraints related directly to stage weight (e.g., TWR and axial load constraints) is restricted in the conventional sizing algorithms, the unified algorithm allows for flexible manipulation of constraints. Therefore, the required constraints in Eq. (6) are properly imposed in the unified algorithm but not considered in the conventional sizing method. Table 3 and Fig. 5 present the case study results for Launcher 2.

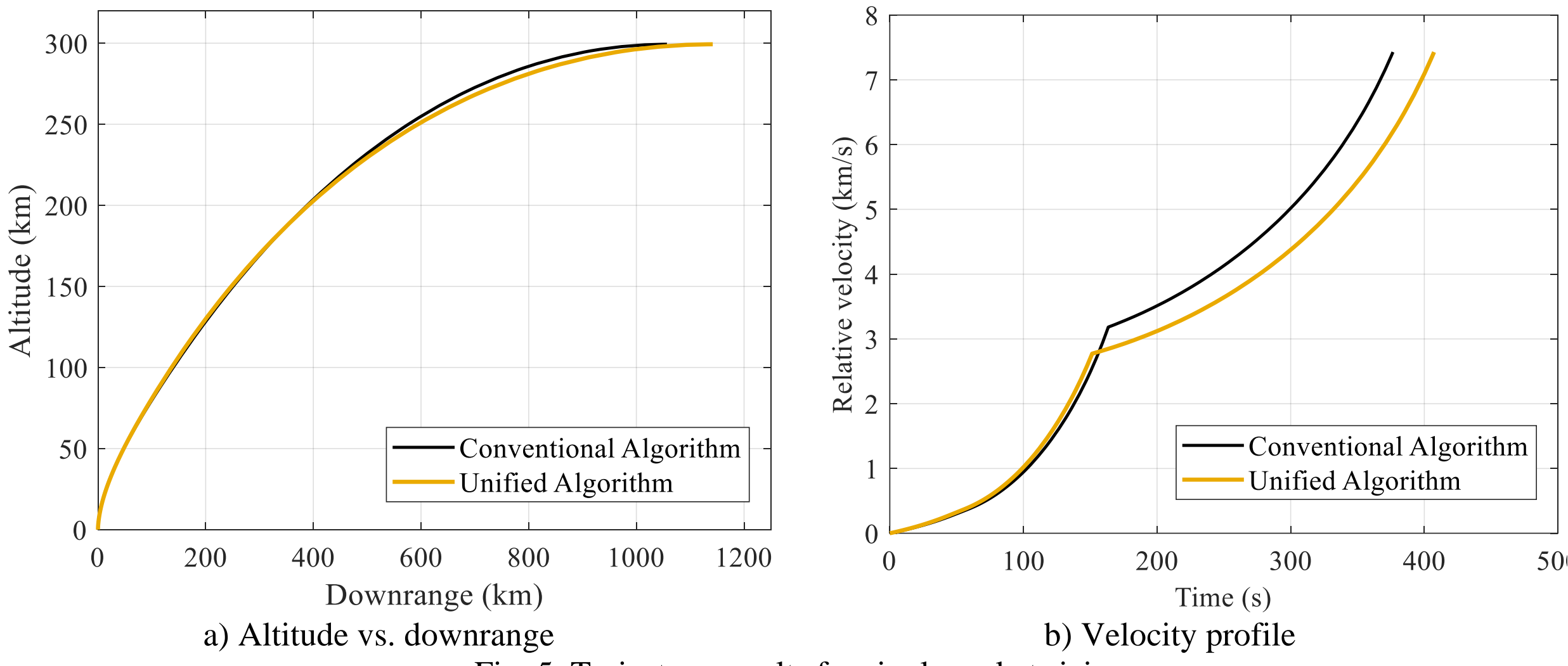

a) Altitude vs. downrange　　b) Velocity profile

Fig. 5. Trajectory results for single rocket sizing

The unified method finds better results, improving the objective function by 3.7% compared to the conventional method in a faster and more tractable manner. The conventional method converges within 5 iterations (the iteration ends when the velocity loss differences between sizing and trajectory modules become smaller than 0.05%), where the average computational time for each iteration is 23 seconds. On the other hand, the unified method provides the optimal solution in a single execution, taking 43.6 seconds. Note that the number of iterations and computation time may vary depending on the initial guesses.

Table 3. Single rocket sizing results

| | $m_{s_1}$, kg | $m_{s_2}$, kg | $m_L$, kg | $m_0$, kg | $m_0/m_L$ |
|---|---|---|---|---|---|
| Unified Algorithm | 4,293 | 1,419 | 5,766 | 157,458 | 27.01 |
| Conventional Algorithm | 4,635 | 1,181 | 5,766 | 161,719 | 28.05 |

## B. Case study 2: rocket family

The case study on a rocket family consists of five scenarios depending on the total number of launches, number of certified reuses ($n_r$), and cost factor for upgrade ($f_u$). We analyze the effect of total launches for each launcher through Cases 1-3 and conduct the sensitivity analysis of $n_r$ and $f_u$ by comparing Cases 2, 4-A, and 4-B. We assume that the reusable stage of Launcher 1 is designed to be reused 15 times in Cases 1-3, while it is upgraded to perform 50 reuses in Case 4. The complementary costs (e.g., landing test, structural reinforcement, and refurbishment costs) must be considered when upgrading the reusable stage for an increased $n_r$. However, we ignore these costs in Case 4-A and, alternatively, add a 20% cost increment in Case 4-B. The upgraded costs ( $C^{(1_u)}$ ) are calculated by multiplying the cost factor for the upgrade by the original cost ( $f_u C^{(1)}$ ).

This paper determines the weighting factors in (34) by the number of launches ratio between two launchers as follows.

$$w^{(k)} = \frac{n_l^{(k)}}{\sum_{k=1}^{2} n_l^{(k)}}. \tag{43}$$

Moreover, we apply a marginal cost of 35% and a production-sharing factor of 0.85.

Tables 4 and 5 provide the cost coefficients used in Eq. (31) and the learning factors applied in the case study for the two launchers. All values are determined with million USD and kg based on empirical data and publicly disclosed information [22, 24, 30-32]. In the case of the RLV, we first estimated cost coefficients by presuming $n_r$ of the current version of Falcon 9 is 15 (Cases 1-3, and 4-A) and updated the coefficients for the upgraded version with 50 certified reuses (Case 4-B). The values for stage 2 of launcher 1 are identical to those for stage 1 of launcher 2 as the stages are compatible. Additionally, we assume that the specifications of the engines and fairing mass are given; hence, the costs for those elements are constant. Technically, the fairing mass varies with the payload mass, but it is inefficient in terms of production cost and time to manufacture a fairing for each individual payload. Therefore, a standardized fairing that can accommodate the heaviest payload among the expected demands is typically designed. It is obvious that the cost of a fuselage is directly attributed to the structure mass of the corresponding stage. Even though the engine mass charges a portion of the structure mass, it is a constant value and thus does not affect the cost estimation of the fuselage. Besides, the operation cost is related to GLOW as it covers the entire launch system's operation, whereas the cost for reuse only applies to the first (reusable) stage.

Table 4. Cost coefficients

| Elements | | $m_{ref}$ | Launcher 1 (RLV) | | Launcher 2 (ELV) | |
|---|---|---|---|---|---|---|
| | | | $k$ | $b$ | $k$ | $b$ |
| Stage 1 fuselage | Dev. | $m_{s_1}$ | 2.685 | 0.50 | Same as stage 2 of Launcher 1 | |
| | Pro. | | 0.035 | 0.65 | | |
| Stage 2 fuselage | Dev. | $m_{s_2}$ | 0.731 | 0.55 | 1.023 | 0.55 |
| | Pro. | | 0.0401 | 0.65 | 0.043 | 0.60 |
| Operation | | $m_0$ | 0.0008 | 0.67 | 0.0004 | 0.67 |
| Reuse | | $m_{s_1}$ | 0.0010 | 0.70 | - | - |

In Launcher 1, the learning factor for the cost incurred to reuse a stage is added. We presume the cost increment from repetitive inspection and maintenance for accumulated fatigue is more significant than the cost reduction from the gained experiences, and thus the value is greater than 1. Other than this, traditional values are used. It is important to note that, in an actual design process, these factors should be seriously determined since the effect on the cost is significant with the increasing number of launches.

Table 5. Learning factors

| Elements | Launcher 1 | Launcher 2 |
|---|---|---|
| Stage 1 structure | 0.95 | 0.95 |
| Stage 1 engine | 0.95 | 0.95 |
| Stage 2 structure | 0.95 | 0.95 |
| Stage 2 engine | 0.95 | 0.95 |
| Fairing | 0.95 | 0.95 |
| Operation | 0.90 | 0.90 |
| Reuse | 1.10 | - |

Fig. 6 illustrates the feasible regions, commonality range, independently optimized configurations, and optimal solutions obtained utilizing the proposed procedure for Cases 1-3. The feasible regions are obtained by repeatedly solving the optimal control problem that minimizes and maximizes remaining stage mass for a given compatible stage within the bounds of comparable stages provided by the commonality check. The depicted regions serve to enhance the comprehension of the results but are not necessarily required in the proposed procedure. The sharp boundaries of the feasible regions stem from weight-related constraints, such as the TWR and axial load limits.

The commonality check is performed using the given specifications of the launchers, and as a result, the commonality range is determined by the minimum mass of Launcher 2 and the maximum mass of Launcher 1. Within the range, the simultaneous sizing process finds the optimal configurations for a rocket family that minimizes the specific cost. Furthermore, the results are compared with the independently optimized configurations obtained from the integrated sizing method introduced in Subsection II-B. The detailed values are represented in Table 6.

The $P_{main}$ converges within 14 iterations for all cases (the iteration ends when the difference between two search points in the golden-section search becomes smaller than 0.01%). The computation time for solving $P_{sub}$ varies depending on the launcher type, initial guesses, and the bounds set for the variables. We made rough initial guesses and used the values for every iteration. Solving the optimal control problems for ELV cases takes 20-40 seconds, and the convergence issue rarely happens. In contrast, the computation time for RLV cases ranges from 90 to 290 seconds over 14 iterations, where we occasionally encounter convergence issues and have to manage them by adjusting the bounds correspondingly.

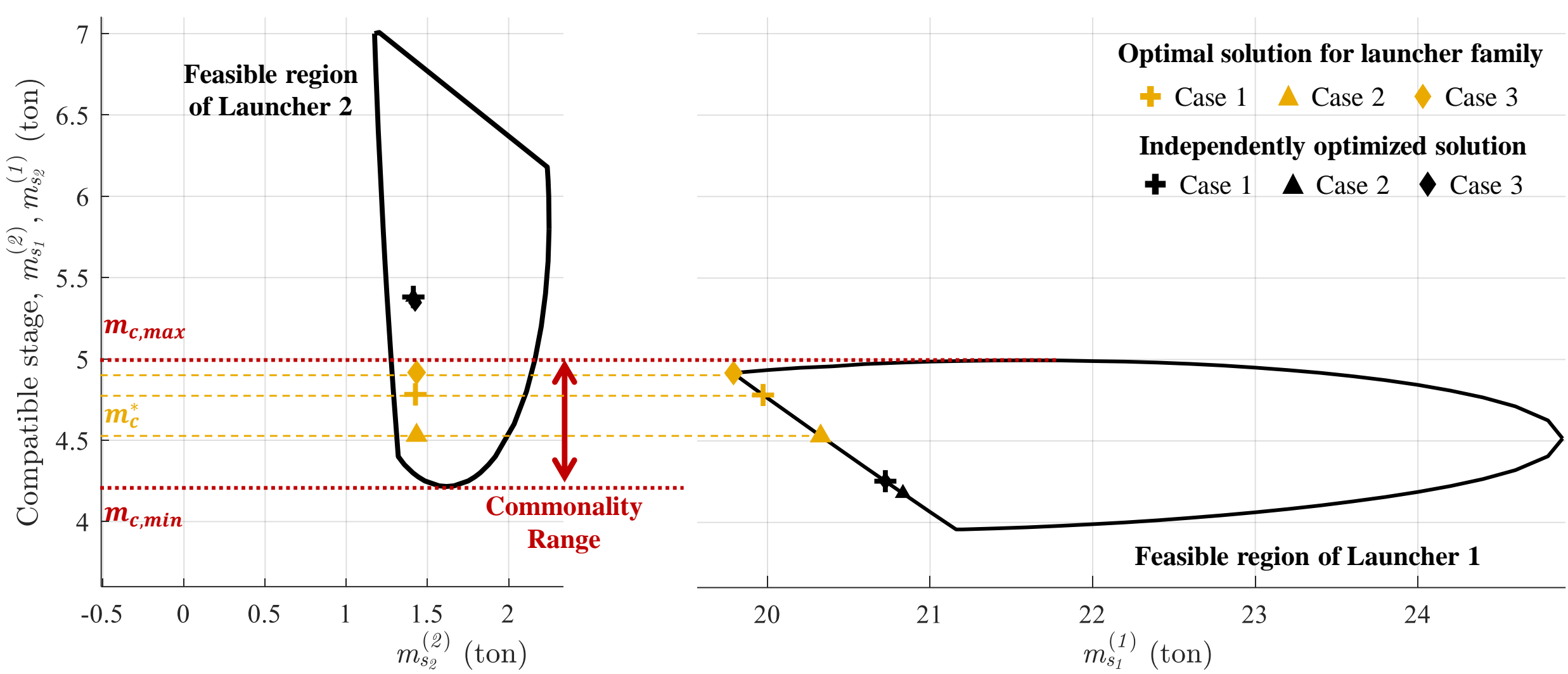

Fig. 6. Feasible regions, commonality range, and optimal configurations of the two launchers

Table 6. Rocket Family sizing results

|  |  | Launcher Index | $m_{s_1}$, kg | $m_{s_2}$, kg | $m_L$, kg | $m_0$, kg | Specific cost, USD/kg | Cost savings, % |
|---|---|---|---|---|---|---|---|---|
| Case 1 $n_l$ = [50, 50], $n_r$ = 15 | Independently Optimized | 1 | 20,725 | 4,253 | 15,858 | 540,520 | 3,167 | - |
| | | 2 | 5,381 | 1,413 | 6,449 | 188,347 | 4,872 | - |
| | Proposed | 1 | 19,973 | 4,782 | 15,769 | 540,520 | 3,008 | 5.0 |
| | | 2 | 4,782 | 1,426 | 6,140 | 171,644 | 4,351 | 10.6 |
| Case 2 $n_l$ = [150, 50], $n_r$ = 15 | Independently Optimized | 1 | 20,833 | 4,179 | 15,819 | 540,520 | 2,231 | - |
| | | 2 | 5,381 | 1,413 | 6,449 | 188,347 | 4,872 | - |
| | Proposed | 1 | 20,327 | 4,530 | 15,889 | 540,520 | 2,124 | 4.9 |
| | | 2 | 4,530 | 1,433 | 5,967 | 164,594 | 4,237 | 12.8 |
| Case 3 $n_l$ = [50, 150], $n_r$ = 15 | Independently Optimized | 1 | 20,725 | 4,253 | 15,858 | 540,520 | 3,167 | - |
| | | 2 | 5,346 | 1,424 | 6,440 | 187,598 | 3,976 | - |
| | Proposed | 1 | 19,791 | 4,917 | 15,628 | 540,589 | 2,971 | 6.2 |
| | | 2 | 4,917 | 1,433 | 6,228 | 175,643 | 3,654 | 8.0 |
| Case 4-A $n_l$ = [150, 50], $n_r$ = 50, $f_u$ = 1.0 | Independently Optimized | 1 | 20,868 | 4,154 | 15,804 | 540,509 | 2,120 | - |
| | | 2 | 5,381 | 1,413 | 6,449 | 188,347 | 4,872 | - |
| | Proposed | 1 | 20,327 | 4,530 | 15,889 | 540,520 | 2,014 | 5.0 |
| | | 2 | 4,530 | 1,433 | 5,967 | 164,594 | 4,237 | 12.8 |
| Case 4-B $n_l$ = [150, 50], $n_r$ = 50, $f_u$ = 1.2 | Independently Optimized | 1 | 20,831 | 4,180 | 15,820 | 540,529 | 2,229 | - |
| | | 2 | 5,381 | 1,413 | 6,449 | 188,347 | 4,872 | - |
| | Proposed | 1 | 20,327 | 4,530 | 15,889 | 540,520 | 2,121 | 4.8 |
| | | 2 | 4,530 | 1,433 | 5,967 | 164,594 | 4,237 | 12.8 |

In all cases, the payload weight exceeds the requirement. In the case of the RLV, the increment is minor because the payload ratio, a fraction of the payload mass in the GLOW, is relatively small, making it inefficient to increase the GLOW. For this reason, the GLOW is constrained by the maximum TWR in all cases. Thus, the GLOW remains almost constant, with only the mass distribution varying, regardless of the method and case. The compatible stage of launcher 2 ( $m_{s_1}^{(2)}$ ) becomes considerably lighter than the independently optimized mass for all cases, while only slight variations occur in the remaining stage ( $m_{s_2}^{(2)}$ ), reducing the GLOW and payload mass. On the other hand, different aspects are observed in the RLV cases: the remarkable increase and decrease in both the compatible stage ( $m_{s_2}^{(1)}$ ) and remaining stage ( $m_{s_1}^{(2)}$ ), respectively, take place. Consequently, the configurations adequately designed for a rocket family can achieve cost savings in all cases.

Fig. 6 and Table 6 show that the independently optimized configuration for the minimum specific cost differs by three factors ($n_l$, $n_r$, and $f_u$). The cost awareness in the ELV design leads to a minor effect on optimal design over $n_r$ change (differences of only a few kilograms in the stage/payload mass among the three cases). However, noticeable differences are observed in the RLV, where the sizing results vary with $n_l$, particularly in the region of low launch counts [22]. This discrepancy comes from the reusability. The costs for both stages of an ELV decrease similarly along the increasing $n_l$, whereas, in an RLV, the cost reduction in the reusable stage is greater than in the other stage as the production cost for the reusable stage is nonrecurring. Namely, $n_l$ and $n_r$ play a pivotal role in an RLV design. Fig. 7 illustrates the objective function versus compatible stage mass within the commonality range for all cases, clearly demonstrating the convex nature of the objective function and the impact of the combination on its shape.

Fig. 6 shows that as the portion of $n_l$ in the particular launcher increases, the optimal solution becomes closer to the independently optimized one. This indicates that the combination of $n_l$ substantially affects the optimal mass configurations in a rocket family since the costs for a compatible stage are distributed into each launcher in proportion to $n_l$ according to Eq. (43). It should be emphasized that the distribution rule, Eq. (43), could be replaced following the design principles. It is obvious that the increased $n_l$ effectively lowers the specific cost in both approaches regardless of the effect of $n_l$ on the optimal mass configuration. As the production quantity increases, the nonrecurring costs (e.g., development cost) are spread out, and the learning effect intensifies, reducing the costs per unit.

Meanwhile, upon comparing Cases 2, 4-A, and 4-B, we can notice that $n_r$ has no discernible influence on the sizing results for a rocket family because only the costs associated with reuses alter, which is not a determining factor in optimizing the configuration of a rocket family. Concerning the specific cost, we obtain the results in Case 4-B similar to Case 2, when we set the complementary cost for the upgrade to be 20%. Excluding this increment, as shown in Case 4-A, results in a cost reduction of approximately 5.4% in the RLV, which is the best case but unrealistic. From this analysis, we can reasonably estimate that the complementary cost for the upgrade ranges from 0% to 20%. It can be reinterpreted that if the upgrade incurs additional costs of more than 20%, increasing the certified number of reuses is ineffective in terms of cost.

Overall, we can conclude that specific costs in a rocket family, as well as the optimal configuration, are highly dependent on the number of launches, emphasizing the importance of precisely predicting the total number of launches in a rocket family design.

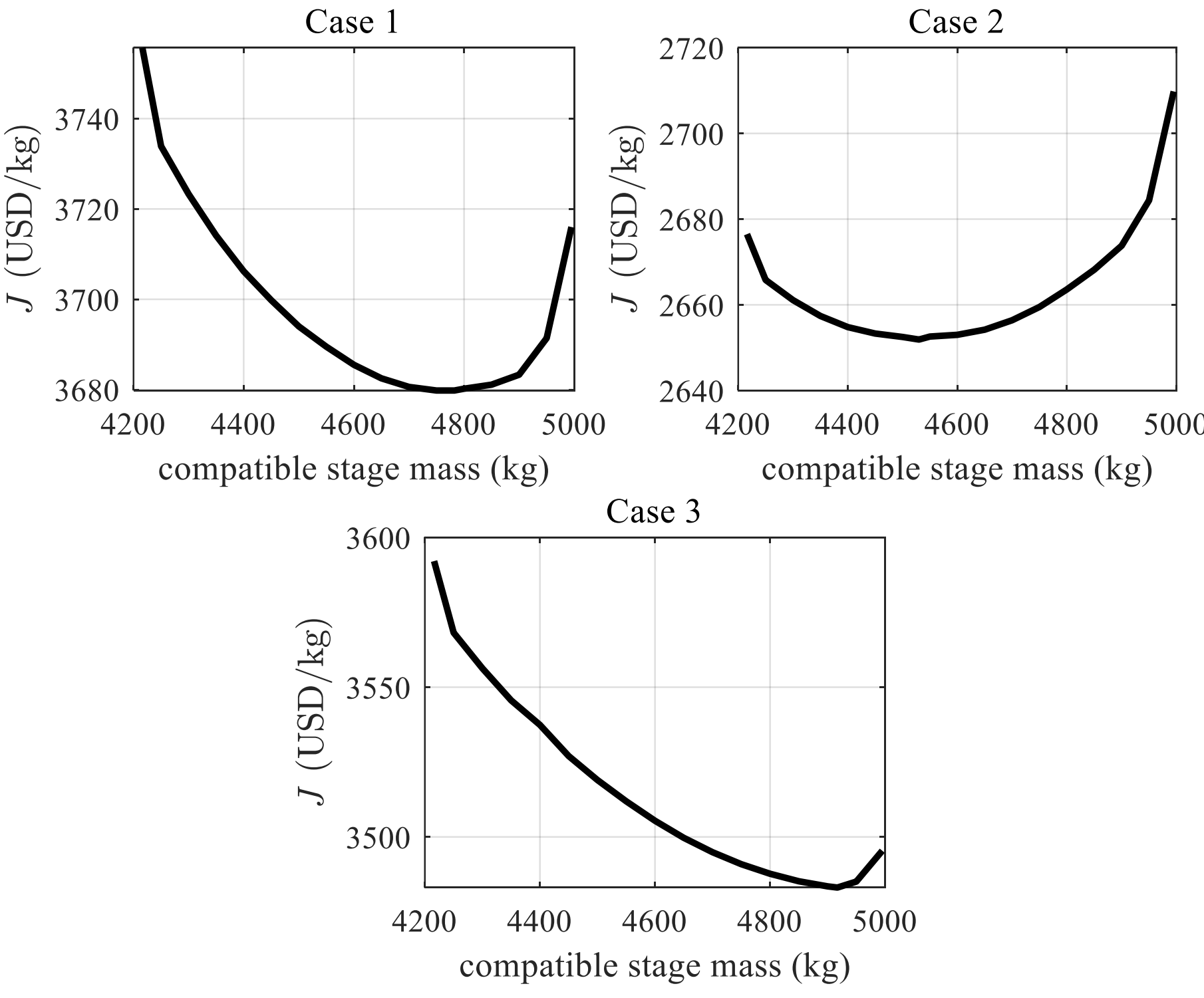

Fig. 7. Objective function versus compatible stage mass for Cases 1-3

The proposed procedure also provides trajectories as outcomes, as shown in Figs. 8-9, in which the independently optimized trajectories are included. Despite the identical target orbit, it is unreasonable to directly compare the two trajectories due to the difference in payload mass. Nevertheless, several noteworthy aspects can be found. Considering that the payload and GLOW are similar to those independently optimized in the case of RLV, we can conclude that the mass configuration evidently influences the entire trajectory. However, the trajectories in Phase 1 are almost identical, with only a different burn time, indicating that the mass of the RLV's first stage may not be the dominant factor in Phase 1. It indicates that different mass configurations with the same GLOW result in varying trajectories, while only GLOW plays a key role in Phase 1 for an RLV. Conversely, only a negligible discrepancy is observed in the two paths of the ELV, even though all mass values are changed, which is a somewhat unexpected result given the substantial mismatch in velocity profiles. This is because the differences in flight path angle and velocity affect the altitude and downrange equally, resulting in similar paths accidentally.

In this case study, the optimal compatible stage mass lies between the independently optimized values for each launcher. As a result, the first stage in both launchers becomes lighter, enforcing the second stage accounting for a larger portion of velocity increment for the orbit insertion. It is observed that the stage configuration affects the engine thrust sequence and, accordingly, the flight time, particularly in an RLV. The burn time for the ELV is directly determined by the propellant mass, which is straightforward. In constraints, the thrust sequence during the descent phases of the RLV is difficult to predict intuitively only with the mass configuration. The total impulse required for the descent phase is mutually dependent on the ascent trajectory and the constraints, necessitating trajectory optimization for accurate sizing.

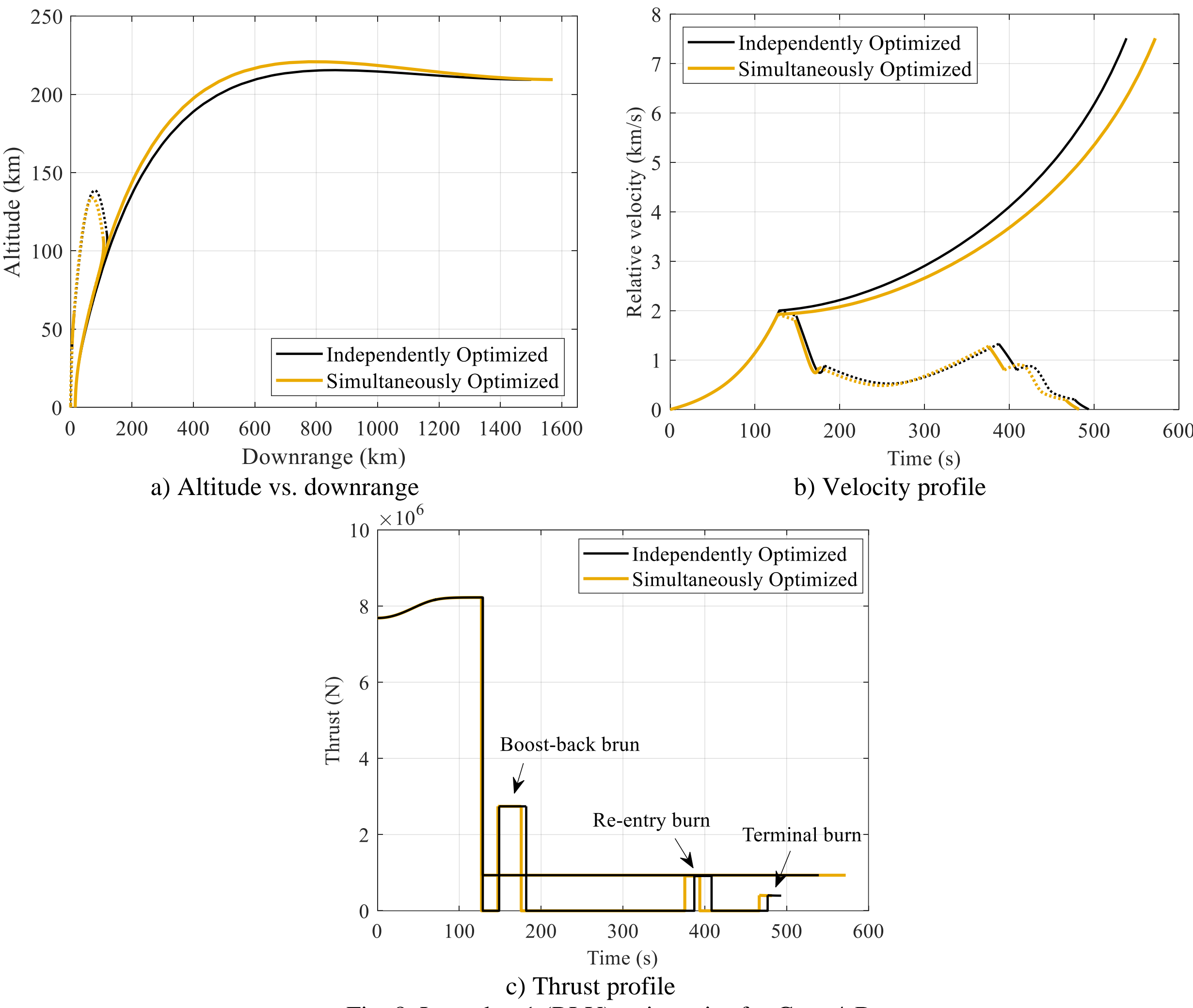

a) Altitude vs. downrange

b) Velocity profile

c) Thrust profile

Fig. 8. Launcher 1 (RLV) trajectories for Case 4-B

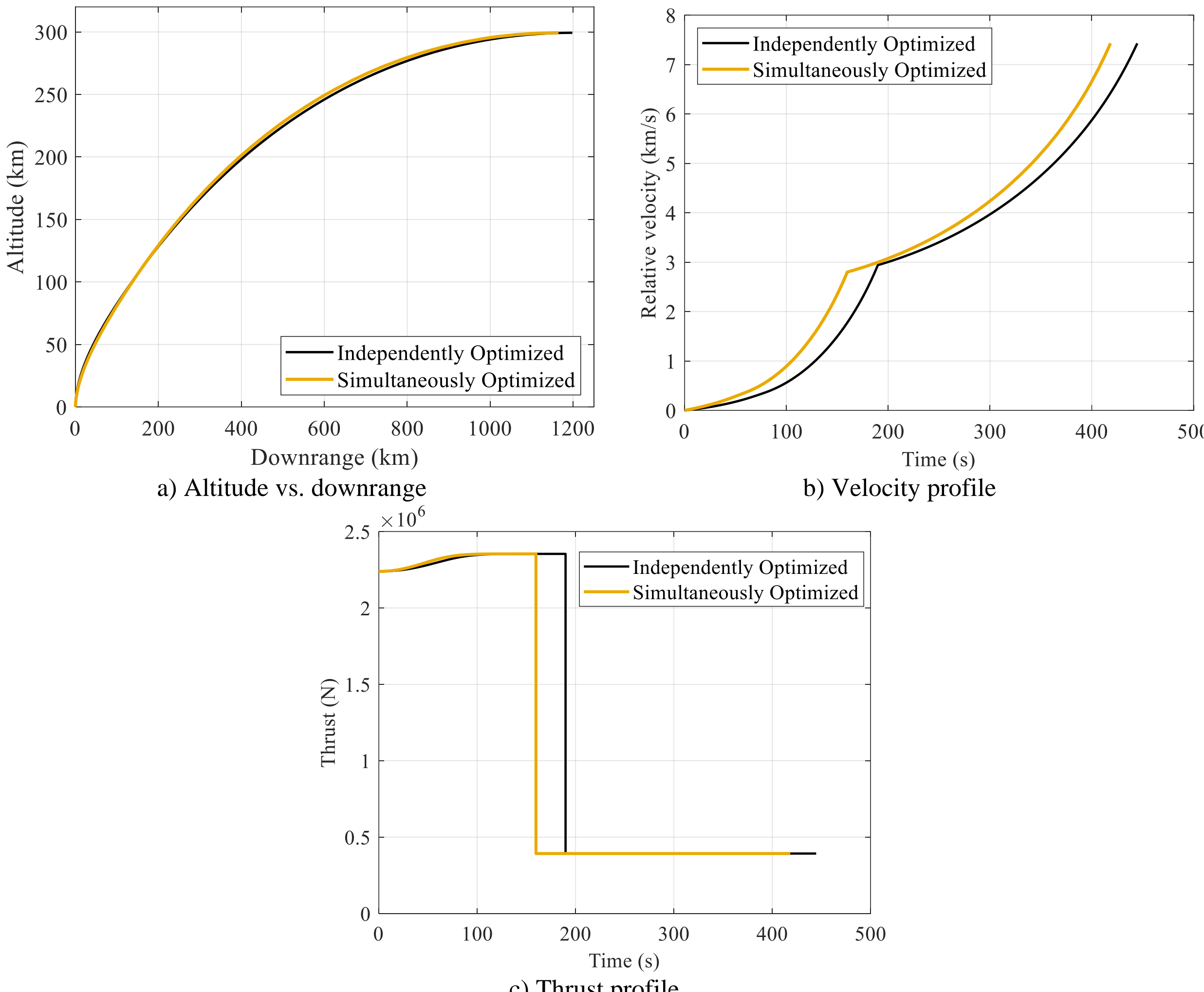

a) Altitude vs. downrange

b) Velocity profile

c) Thrust profile

Fig. 9. Launcher 2 (ELV) trajectories for Case 4-B

## V. Conclusion

This paper proposes a sizing procedure for a rocket family consisting of two launchers in which payload capability differs by considering the commonalities between the launchers. Instead of individually or sequentially designing each launcher, this procedure integrates all launchers into one design process to leverage the full advantage of sharing compatible stages, leading to lower development and manufacturing costs for the compatible stages. Unlike the conventional sizing algorithms for a single rocket, whose iterative loops are necessary and the optimality gap is relatively large, the proposed procedure transforms the staging problem into a single optimal control problem by imposing additional constraints for stage properties. Based on this method, a commonality check is performed to determine the commonality range by obtaining the min/max stage weight of a compatible stage of each launcher. Meanwhile, a cost model for a rocket family is established by considering commonality and referring to a cost model for an independent rocket. Afterward, an optimal control problem that minimizes the specific launch cost of the payload (USD/kg) for a given compatible stage weight is formulated for each launcher. As a result, the problem is transformed into an optimization problem with single-dimensional search spaces over the compatible stage weight. Lastly, an optimal compatible stage weight is calculated via a direct search method.

A case study on a two-member rocket family containing the preliminarily designed ELV and the existing model for RLV, Falcon 9, is carried out, demonstrating the effectiveness of the proposed procedure. It is shown that the RLV and ELV can save costs of approximately 5~13% depending on the number of launches and reuses by sharing the compatible stage and finding optimal configurations. Mass configurations alter considerably in most cases compared to the independently optimized solutions, resulting in lightening the payload weight and GLOW for the ELV, whereas the optimized payload weight and GLOW of the RLV remain similar. Furthermore, the proposed procedure optimizes not only the mass configurations but also the trajectories for all family members, which is another core contribution of this paper. The potential future works include a comprehensive design for a rocket family from the multidisciplinary design optimization perspective, a simultaneous sizing of a rocket family consisting of more than two launchers, or consideration of notable types of launch systems, such as the air launch and spin launch. Moreover, addressing practical issues and improving accuracy would be interesting future work to enhance the algorithm’s fidelity.

## Acknowledgments

The authors would like to thank Nicholas Gollins and Erika Kim for their editorial review and thoughtful suggestions for improvement.